\documentclass[11pt]{article}
\usepackage{amsmath}
\usepackage{amsfonts}
\usepackage[bookmarksnumbered=true]{hyperref}

\newtheorem{theorem}{Theorem}[section]
\newtheorem{lemma}[theorem]{Lemma}
\newtheorem{proposition}[theorem]{Proposition}

\newtheorem{definition}[theorem]{Definition}
\newtheorem{remark}[theorem]{Remark}

\begin{document}

\title{Forward-Backward Doubly Stochastic Differential Equations
with Random Jumps and Stochastic Partial Differential-Integral
Equations\thanks{This work is supported by National Natural Science Foundation of China Grant
10771122, Natural Science Foundation of Shandong Province of China
Grant Y2006A08 and National Basic Research Program of China (973
Program, No.2007CB814900)}}
\author{Qingfeng Zhu$^{\rm a}$ and Yufeng Shi$^{\rm b}$\thanks{Corresponding author, E-mail: yfshi@sdu.edu.cn}\\
{\small $^{\rm a}$ School of Statistics and Mathematics, Shandong
University of Finance},\\
{\small Jinan 250014, China}\\
{\small$^{\rm b}$School of Mathematics, Shandong University, Jinan 250100, China}}
\maketitle

\begin{abstract}In this paper, we study
forward-backward doubly stochastic differential equations
driven by Brownian motions and Poisson process (FBDSDEP in short).
Both the probabilistic interpretation for the solutions to
a class of quasilinear stochastic partial differential-integral
equations (SPDIEs in short) and stochastic Hamiltonian systems
arising in stochastic optimal control problems with random jumps are
treated with FBDSDEP. Under some monotonicity assumptions, the
existence and uniqueness results for measurable solutions of FBDSDEP
are established via a method of continuation. Furthermore, the
continuity and differentiability of the solutions of FBDSDEP
depending on parameters is discussed. Finally, the probabilistic
interpretation for the solutions to a class of quasilinear SPDIEs is
given.\\
\indent{\it keywords:} forward-backward doubly stochastic differential equations,\\
stochastic partial differential-integral equations, random measure, Poisson process
\end{abstract}

\section{Introduction}\label{sec:1}

Nonlinear backward stochastic differential equations with Brownian
motion as noise sources (BSDEs in short) were first introduced by
Pardoux and Peng \cite{PP1}. By virtue of BSDEs, Peng \cite{P} has given a
probabilistic interpretation (nonlinear Feynman-Kac formula) for the
solutions of semilinear parabolic partial differential equations
(PDEs in short), for more detailed, the reader is referred to
Darling and Pardoux \cite{DP}, Pardoux and Zhang \cite{PZ} and so on. Fully
coupled forward-backward stochastic differential equations (FBSDEs
in short) can provide a probabilistic interpretation for the
solutions to a class of quasilinear parabolic and elliptic PDEs (cf
Pardoux and Tang \cite{PT}) and have been investigated deeply. FBSDEs
were studied first by Antonelli \cite{A} and Ma et al. \cite{MPY}, Cvitanic and
Ma \cite{CM} and used to hedge options involved in a large investor in
financial market. Consequently, FBSDEs were developed in Hu and Peng
\cite{HP}, Peng and Wu \cite{PW} and Peng and Shi \cite{PS1} and so on.

A class of backward doubly stochastic differential equations (BDSDEs
in short) was introduced by Pardoux and Peng \cite{PP2} in 1994, in order
to provide a probabilistic interpretation for the solutions to a
class of quasilinear stochastic partial differential equations
(SPDEs in short). Due to their important significance to SPDEs , the
researches for BDSDEs have been in the ascendant (cf. Bally and
Matoussi \cite{BM}, Zhang and Zhao \cite{ZZ}, Ren et al. \cite{RLH}, Hu and Ren \cite{HR}
and their references). In 2003, Peng and Shi \cite{PS2} have introduced a
type of time-symmetric forward-backward stochastic differential
equations, which is a generalization of stochastic Hamilton system.
Recently Zhu et al. \cite{ZSG} have extended the results in \cite{PS2} to
different dimensional forward-backward doubly stochastic
differential equations (FBDSDEs in short) and weakened the monotone
assumptions. However, the theory of FBDSDEs has not been
investigated enough so far.

BSDEs driven by Brownian motions and Poisson process (BSDEP in
short) was first discussed by Tang and Li \cite{TL}. After then Situ \cite{S}
obtained an existence and uniqueness result with non-Lipschitz
coefficients for BSDEP, so as to get the probabilistic
interpretation for solutions of  partial differential-integral
equations (PDIEs in short). Using this kind of BSDEP Barles et al.
\cite{BBP} and Yin and Mao \cite{YM} proved that there exists a unique viscosity
solution for a system of parabolic integral-partial differential
equations. Fully coupled forward-backward stochastic differential
equations with Poisson process (FBSDEP in short) was discussed by Wu
\cite{W} and Yin and Situ \cite{YS}. Then FBSDEP were used to study the
linear quadratic optimal control problems with random jump by Wu and
Wang \cite{WW} and the nonzero-sum differential games with random jumps
by Wu and Yu \cite{WY}. Recently BDSDEs driven by Brownian motions and
Poisson process (BDSDEP in short) was discussed by Sun and Lu \cite{SL}.

Because of their important significance to SPDEs, it is necessary to
give intensive investigation to the theory of FBDSDEs. In this paper
we study FBDSDEs driven by Brownian motions and Poisson process
(FBDSDEP in short), which generalize the so-called time-symmetric
forward-backward stochastic differential equations introduced by
Peng and Shi \cite{PS2}. FBDSDEP can provide more extensive frameworks for
the probabilistic interpretation (nonlinear stochastic Feynman-Kac
formula) for the solutions to a class of quasilinear stochastic
partial differential-integral equations (SPDIEs in short) and
stochastic Hamiltonian systems arising in stochastic optimal control
problems with random jumps. Under some monotonicity assumptions, we
establish the existence and uniqueness results for measurable
solutions of FBDSDEP by means of a method of continuation
systemically introduced by Yong \cite{Y}. Then we discuss the continuity
and differentiability of the solutions to FBDSDEP depending on
parameters. Furthermore, by virtue of FBDSDEP, we give the
probabilistic interpretation for the solutions to a class of
quasilinear SPDIEs. Finally, we discuss a doubly stochastic
Hamiltonian system.

The paper is organized as follows. In Section \ref{sec:2}, some preliminary
results are given. Section \ref{sec:3} is devoted to proving the existence and
uniqueness result for FBDSDEP. In Section \ref{sec:4}, the continuity and
differentiability of the solutions to FBDSDEP depending on
parameters is discussed. In Section \ref{sec:5}, the probabilistic
interpretation for the solutions to a class of quasilinear SPDIEs is
given by virtue of this class of FBDSDEP. Finally, in Section \ref{sec:6}, the
above results are applied to a doubly stochastic Hamiltonian system.

\section{Preliminary}\label{sec:2}

Let $(\Omega,{\cal F},P) $ be a complete probability space, and
$[0,T]$ be a fixed arbitrarily large time duration throughout this
paper. We suppose $\{{\cal F}_t\}_{t \geq 0}$ is generated by the
following three mutually independent processes:

(i)\hspace{0.1cm}Let $\left\{ W_t;0\leq t\leq T\right\}$ and
$\left\{ B_t;0\leq t\leq T\right\} $ be two standard Brownian
motions defined on $( \Omega ,{\cal F},P) $, with values
respectively in $\mathbb{R}^d$ and in $\mathbb{R}^l$.

(ii)\hspace{0.1cm}Let $N$ be a Poisson random measure, on
$\mathbb{R}_{+} \times Z$, where $Z \subset \mathbb{R}^r$ is a
nonempty open set equipped with its Borel field ${\cal B}(Z)$, with
compensator $\widehat{N}(dzdt)=\lambda(dz)dt$, such that $\widetilde{N}(A \times [0,t])=(N-\widehat{N})(A
\times [0,t])_{t\geq 0}$ is a martingale for all $A \in {\cal B}(Z)$
satisfying $\lambda(A) < \infty$. $\lambda$ is assumed to be a
$\sigma$-finite measure on $(Z,{\cal B}(Z))$ and is called the
characteristic measure.

Let ${\cal N}$ denote the class of $P$-null elements of ${\cal F}$.
For each $t\in \left[ 0,T\right] $, we define $ {\cal F}_t\doteq
{\cal F}_t^W\vee {\cal F}_{t,T}^B\vee {\cal F}_t^N,$ where for any
process $\{\eta_t\},{\cal F}_{s,t}^\eta=\sigma
\left\{\eta_r-\eta_s;s\leq r\leq t\right\}\vee {\cal N}, {\cal
F}_t^\eta={\cal F}_{0,t}^\eta $. Note that the collection $\left\{
{\cal F}_t,t\in \left[ 0,T\right] \right\} $ is neither increasing
nor decreasing, and it does not constitute a classical filtration.

We introduce the following notations:
\begin{eqnarray*}
M^2(0,T;\mathbb{R}^n)&=&\{v_t,0\leq t\leq T,\ \mbox {is an}\
\mathbb{R}^n\mbox{-valued},\ {\cal F}_t\mbox{-measurable process}\\
&&\mbox{such that} \ E\int_{0}^T|v_{t}|^{2}dt<\infty\},\\
F^2_N(0,T;\mathbb{R}^n)&=&\{k_t,0\leq t\leq T,\ \mbox {is an}\
\mathbb{R}^n\mbox{-valued},\ {\cal F}_t\mbox{-measurable process}\\
&&\mbox{such
that} \ E\int_{0}^T\int_{Z}|k_{t}(z)|^{2}\lambda(dz)dt<\infty\},\\
L^2_{\lambda(\cdot)}(\mathbb{R}^n)&=&\{k(z),\ k(z)\ \mbox{is an}\
\mathbb{R}^n\mbox{-valued}, \ {\cal B}(Z)\mbox{-measurable function}\\
&&\mbox{such that}\ \|k\|=(\int_Z|k(z)|^2\lambda(dz))^{1/2}<\infty\},\\
L^{2}(\Omega, {\cal F}_T,P;\mathbb{R}^{n})&=&\{\xi,\ \xi\ \mbox{is
an}\ \mathbb{R}^n\mbox{-valued},\  {\cal F}_T\mbox{-measurable
random }\\
&&\mbox{ variable such that}\ E|\xi|^2 <\infty\}.
\end{eqnarray*}
For a given $u\in M^2( 0,T;\mathbb{R}^d) $ and $v\in M^2(
0,T;\mathbb{R}^l) $, one can define the (standard) forward It\^o's
integral $\int_0^{\cdot }u_sdW_s$ and the backward It\^o's
integral $\int_{\cdot }^Tv_sdB_s$. They are both in $M^2(
0,T;\mathbb{R})$ (see \cite{PP2}). We use the usual inner product $\langle
\cdot ,\cdot \rangle $ and Euclidean norm $ | \cdot | $ in
$\mathbb{R}^n$, $\mathbb{R}^m$, $\mathbb{R}^{m\times l}$ and
$\mathbb{R}^{n\times d}.$ All the equalities and inequalities
mentioned in this paper are in the sense of $dt\times dP$ almost surely on $\left[ 0,T\right] \times \Omega $.

Consider the following BDSDE with Brownian motions and Poisson
Process (BDSDEP in short):
\begin{equation}\label{eq:1}
\left\{
\begin{array}{lll}
dP_t&=&F(t,P_t,Q_t,K_t)dt+G( t,P_t,Q_t,K_t)dB_t\\
&&+Q_tdW_t+\displaystyle\int_ZK_{t_-}(z)\widetilde{N}(dzdt),\\
P_T&=&\xi,
\end{array}\right.
\end{equation}
where
\begin{eqnarray*}
&&F:\Omega \times [0,T] \times \mathbb{R}^{m} \times \mathbb{R}^{m
\times d} \times L^2_{\lambda(\cdot)}(\mathbb{R}^m)\rightarrow
\mathbb{R}^{m},\\
&&G:\Omega \times [0,T] \times
\mathbb{R}^{m} \times \mathbb{R}^{m \times d} \times
L^2_{\lambda(\cdot)}(\mathbb{R}^m)\rightarrow \mathbb{R}^{m \times
l}.
\end{eqnarray*}

%\newdefinition
\begin{definition}\label{def:2.1}
A triple of  ${\cal F}_t$-measurable stochastic processes $(P,Q,K)$
is called a solution to BDSDEP (\ref{eq:1}) if $(P$, $Q,K)\in M^2([0,T];
\mathbb{R}^{m+m\times d})\times F^2_N(0,T;\mathbb{R}^{m})$ and
satisfies BDSDEP (\ref{eq:1}).
\end{definition}

We assume that

\begin{enumerate}
\item[(H1)]
\begin{enumerate}
\item[(i)] $\xi \in L^{2}( \Omega ,{\cal F}_T,P;\mathbb{R}^{m})$;
\item[(ii)] $F$ is $
{\cal F}_t$-progressively measurable and satisfies
$F(\omega,t,0,0,0)\in M^2(0$, $T;\mathbb{R}^{m})$;
\item[(iii)] $G$ is ${\cal F}_t$-progressively measurable and satisfies
$G(\omega,t,0,0,0)\in M^2(0$, $T;\mathbb{R}^{m\times l})$;
\item[(iv)] $F$ and $G$ satisfy Lipschitz conditions
to $P,Q,K$, that is, there exist $c >0$ and $0<\gamma<1$ such that
\begin{eqnarray*}
&&|F(t,P_{1},Q_{1},K_{1})-F(t,P_{2},Q_{2},K_{2})|^2\\
&\leq& c(|P_1-P_2|^2+|Q_1-Q_2|^2+\|K_1-K_2\|^2),\\
&& |G(t,P_{1},Q_{1},K_{1})-G(t,P_{2},Q_{2},K_{2})|^2\\
&\leq& c|P_1-P_2|^2+\gamma(|Q_1-Q_2|^2+\|K_1-K_2\|^2).
\end{eqnarray*}
\end{enumerate}
\end{enumerate}

In order to attain our results, we give the following Proposition
\ref{pro:2.2} and Proposition \ref{pro:2.3} appeared in \cite{SL}
and Proposition \ref{pro:2.4} appeared in \cite{S}.

\begin{proposition}\label{pro:2.2}
Under the assumption (H1), BDSDEP (\ref{eq:1}) has a unique
solution $(P,Q,K)\in M^2([0,T];\mathbb{R}^{m+m\times d})\times
F^2_N(0,T;\mathbb{R}^{m}).$
\end{proposition}

\begin{proposition}\label{pro:2.3}
Let $\alpha \in M^2([0,T];\mathbb{R}^m),\beta \in
M^2(0,T;\mathbb{R}^{m}),\gamma \in M^2(0,T$;\\ $\mathbb{R}^{m \times l})$,
$\delta \in M^2(0,T;\mathbb{R}^{m \times d})$ be such that:
\begin{eqnarray*}
\alpha_t&=&\alpha_0+\int_0^t\beta_s ds+\int_0^t \gamma_s dB_s+\int_0^t\delta_sdW_s\\
&&+\int_0^t\int_ZK_{s_-}(z)\widetilde{N}(dzds), \
0\leq t\leq T.
\end{eqnarray*}
Then
\begin{eqnarray*}
|\alpha_t|^2&=&|\alpha_0|^2+2\int_0^t\langle\alpha_s,\beta_s\rangle ds
+2\int_0^t\langle\alpha_s,\gamma_sdB_s\rangle\\
&&+2\int_0^t\langle\alpha_s,\delta_sdW_s\rangle
+2\int_0^t\int_Z\langle\alpha_s,K_{s_-}(z)\rangle\widetilde{N}(dzds)\\ &&-\int_0^t|\gamma_s|^2ds+\int_0^t|\delta_s|^2ds+\int_0^t\|K_{s}\|^2ds,
\end{eqnarray*}
\begin{eqnarray*}
E|\alpha_t|^2&=&E|\alpha_0|^2
+2E\int_0^t\langle\alpha_s,\beta_s\rangle ds -E\int_0^t|\gamma_s|^2ds\\
&&+E\int_0^t|\delta_s|^2ds+E\int_0^t\|K_{s}\|^2ds.
\end{eqnarray*}
\end{proposition}

\begin{proposition}\label{pro:2.4}
Let $\alpha \in M^2([0,T];\mathbb{R}^n),\beta \in
M^2(0,T;\mathbb{R}^{n}),\gamma \in M^2(0,T;\mathbb{R}^{n \times l})$,
$\delta \in F_N^2(0,T;\mathbb{R}^{n})$ be such that:
\begin{eqnarray*}
\alpha_t=\alpha_0+\int_0^t\beta_s ds+\int_0^t\gamma_sdW_s+\int_0^t\int_Z\delta_{s_-}(z)\widetilde{N}(dzds),
\ \ \ \ \ 0\leq t\leq T.
\end{eqnarray*}
Then for all $u(s,x)\in C^{1,2}([0,T]\times \mathbb{R}^n;
\mathbb{R}^m)$, one has that
\begin{eqnarray*}
u(t,\alpha_t)-u(0,\alpha_0)&=&\int_0^t{\cal L}u(s,\alpha_s)ds
+\int_0^t \nabla u(s,\alpha_s)\gamma_sdW_s\\
&&+\int_0^t\int_Z(u(s_-,\alpha_{s_-}+\delta(s_-,z))-u(s_-,\alpha_{s_-}))\widetilde{N}(dzds),
\end{eqnarray*}
where
\begin{eqnarray*}
{\cal L}u =\left(\begin{array}{c} Lu_1\\\vdots\\
Lu_m\end{array}\right),
\end{eqnarray*}
with
\begin{eqnarray*}
&&Lu_k(t,x)\\
&=&\frac{\partial u_k}{\partial t}(t,x)
+\sum\limits_{i=1}^n\beta_{i}(t)\frac{\partial u_k}{\partial
x_i}(t,x)+\frac{1}{2}\sum\limits_{i,j=1}^n(\gamma
\gamma*)_{ij}(t)\frac{\partial^2 u_k}{\partial x_i\partial x_j}(t,x)\\
&&+\displaystyle\int_Z(u_k(t,x+\delta(t,x,z))-u_k(t,x)-\sum\limits_{i=1}^n
\delta_i(t,x,z)\displaystyle\frac{\partial u_k}{\partial x_i}(t,x)
)\lambda(dz),\\
 &&k=1,\cdots,m.
\end{eqnarray*}
\end{proposition}

\section{The existence and uniqueness theorem of FBDSDEP}
\label{sec:3}

Consider the following forward-backward doubly stochastic
differential equations with Brownian motions and Poisson process
(FBDSDEP in short):
\begin{equation}\label{eq:2}
\left\{
\begin{array}{l}
dX_t=f(t,X_t,P_t,Y_t,Q_t,K_t)dt+g(
t,X_t,P_t,Y_t,Q_t,K_t)dW_t-Y_tdB_t\\
\qquad\quad+\displaystyle\int_Zh(t_-,X_{t_-},P_{t_-},Y_{t_-},Q_{t_-},K_{t_-}(z),z)
\widetilde{N}(dzdt),  \\
dP_t=F(t,X_t,P_t,Y_t,Q_t,K_t)dt+G(
t,X_t,P_t,Y_t,Q_t,K_t)dB_t+Q_tdW_t\\
\qquad\quad+\displaystyle\int_ZK_{t_-}(z)\widetilde{N}(dzdt),\\
X_0=\Psi (P_0),\quad P_T=\Phi ( X_T).
\end{array}\right.
\end{equation}
where
\begin{center}
$\begin{array}{llll} F & :\Omega \times \left[ 0,T\right] \times
\mathbb{R}^n\times \mathbb{R}^m\times \mathbb{R}^{n\times l}\times
\mathbb{R}^{m\times d} \times L^2_{\lambda(\cdot)}(\mathbb{R}^m)
 & \rightarrow & \mathbb{R}^m, \\
f & :\Omega \times \left[ 0,T\right] \times \mathbb{R}^n\times
\mathbb{R}^m\times \mathbb{R}^{n\times l}\times \mathbb{R}^{m\times
d}\times L^2_{\lambda(\cdot)}(\mathbb{R}^m)  &
\rightarrow & \mathbb{R}^n, \\
G & :\Omega \times \left[ 0,T\right] \times \mathbb{R}^n\times
\mathbb{R}^m\times \mathbb{R}^{n\times l}\times \mathbb{R}^{m\times
d}\times L^2_{\lambda(\cdot)}(\mathbb{R}^m)
  & \rightarrow & \mathbb{R}^{m\times l}, \\
g & :\Omega \times \left[ 0,T\right] \times \mathbb{R}^n\times
\mathbb{R}^m\times \mathbb{R}^{n\times l}\times \mathbb{R}^{m\times
d}\times L^2_{\lambda(\cdot)}(\mathbb{R}^m) & \rightarrow &
\mathbb{R}^{n\times d}, \\
h & :\Omega \times \left[ 0,T\right] \times \mathbb{R}^n\times
\mathbb{R}^m\times \mathbb{R}^{n\times l}\times \mathbb{R}^{m\times
d}\times L^2_{\lambda(\cdot)}(\mathbb{R}^m)\times Z & \rightarrow & \mathbb{R}^{n},\\
\Psi & :\Omega \times \mathbb{R}^m\rightarrow \mathbb{R}^n,\quad
\Phi :\Omega \times \mathbb{R}^n\rightarrow \mathbb{R}^m.
\end{array}$
\end{center}
Given an $m\times n$ full-rank matrix $H$. Let us introduce some
notations
\begin{eqnarray*} &&U =(X,P,Y,Q,K),\quad
A( t,U )=(H^{\rm T}F,Hf,H^{\rm T}G,Hg,Hh)(t,U),\\
&&\langle A,U \rangle=\langle X,H^{\rm T}F\rangle+\langle P,Hf
\rangle +\langle Y, H^{\rm T}G \rangle +\langle
Q,Hg\rangle+\langle\langle K,Hh\rangle\rangle,
\end{eqnarray*}
where
\begin{eqnarray*}
H^{\rm T}G&=&(H^{\rm T}G_1\cdots H^{\rm T}G_l),\quad Hg=(
Hg_{1}\cdots Hg_d),\\
\langle\langle K,Hh\rangle\rangle&=&\int_Z\langle
K_t(z),Hh(t,U_{t},z)\rangle\lambda(dz).
\end{eqnarray*}
\begin{definition}\label{de:3.1}
A quintuple of ${\cal F}_t$-measurable stochastic processes
$(X,P,Y,Q,K)$ is called a solution to FBDSDEP (\ref{eq:2}), if $(X,P,Y$,
$Q,K)\in M^{2}( 0,T$;\\ $\mathbb{R}^{n+m+n\times l+m\times d})\times
F^2_N(0,T;\mathbb{R}^m)$ and satisfies FBDSDEP (\ref{eq:2}).
\end{definition}
The following monotonicity conditions are our main assumptions:
\begin{enumerate}
\item[(H2)]$\forall
U=(X,P,Y,Q,K),\bar{U}=(\bar{X},\bar{P},\bar{Y},\bar{Q},\bar{K})\in
\mathbb{R}^{n+m+n\times l+m\times d}\times L^2_{\lambda(\cdot)}(
\mathbb{R}^m)$, $\forall t\in \left[ 0,T\right]$,
\begin{eqnarray*}
&&\langle A( t,U ) -A( t,\bar{U}) ,U -\bar{U}\rangle\\
&\leq& -\mu _{1}(|
H(X-\bar{X})|^{2}+| H(Y-\bar{Y})|^{2})\\
&& -\mu_{2}(|H^{\rm
T}(P-\bar{P})|^{2}+|H^{T}(Q-\bar{Q})|^{2}+\|H^{\rm
T}(K-\bar{K})\|^{2}).
\end{eqnarray*}
\item[(H3)] $\langle \Psi( P) -\Psi( \bar{P}), H^{\rm T}( P-\bar{P})
\rangle \leq -\beta _{2}| H^{\rm
T}(P-\bar{P})|^{2},\quad \forall P,\bar{P}\in \mathbb{R}^{m}, \\
\langle \Phi (X) -\Phi ( \bar{X}), H(X-\bar{X}) \rangle \geq \beta
_{1}| H(X-\bar{X})| ^{2},\quad \forall X,\bar{X}\in \mathbb{R}^{n}.$
\end{enumerate}
Here, $\mu _{1}$, $\mu _{2}$, $\beta _{1}$ and $\beta _{2}$ are
given nonnegative constants with $\mu _{1}+\mu _{2}>0$, $\beta
_{1}+\beta _{2}>0$, $\mu_{1}+\beta _{2}>0$, $\mu _{2}+\beta _{1}>0.$
Moreover we have $\mu _{1}>0$, $\beta _{1}>0$ (resp. $\mu _{2}>0$,
$\beta _{2}>0$ ) when $m>n$ (resp. $n>m$).

We also assume that
\begin{enumerate}
\item[(H4)]For each $U \in \mathbb{R}^{n+m+n\times l+m\times d}\times
 L^2_{\lambda(\cdot)}(\mathbb{R}^m)$, $A( \cdot,U)$
  is an ${\cal F}_{t}$-measurable vector process defined on
$\left[ 0,T\right]$ with $A(\cdot,0) \in M^{2}( 0,T$;
$\mathbb{R}^{n+m+n\times l+m\times d})\times F^2_N(0,T;\mathbb{R}^m)
$, and for each $X \in \mathbb{R}^{n},\Phi ( X) $ is an ${\cal
F}_{T}$-measurable vector process with $\Phi ( 0) \in L^{2}( \Omega
,{\cal F}_{T},P;\mathbb{R}^{m}) $, and for each $P\in
\mathbb{R}^{m}$, $\Psi ( P) $ is an ${\cal F}_{0}$-measurable vector
process with $\Psi ( 0) \in L^{2}( \Omega,{\cal F}
_{0},P;\mathbb{R}^{n}) .$

\item[(H5)]$f,F,g,G,h,\Psi$ and $\Phi $ satisfy the
Lipschitz conditions:  there exist constants $c>0$ and $0<\gamma <1$
such that
 $\forall U=(X,P,Y,Q,K),
\bar{U}=( t,\bar{X},\bar{P},\bar{Y},\bar{Q},\bar{K})$ $\in
\mathbb{R}^{n+m+n\times l+m\times d}\times L^2_{\lambda(\cdot)}(
\mathbb{R}^m),\forall t\in \left[ 0,T\right]$,
\begin{eqnarray*}
&&|f( t,X,P,Y,Q,K) -f( t,\bar{X},\bar{P},\bar{Y},\bar{Q},\bar{K})|^{2}\\
&\leq& c( |X-\bar{X}| ^{2}+|P-\bar{P}|^{2}+|Y-\bar{Y}|^{2}
+|Q-\bar{Q} |^2+\|K-\bar{K}\|^2),\\
&&|F( t,X,P,Y,Q,K) -F( t,\bar{X},\bar{P},\bar{Y},\bar{Q},\bar{K})|^{2}\\
&\leq& c( |X-\bar{X}| ^{2}+|P-\bar{P}|^{2}+|Y-\bar{Y}|^{2}
+|Q-\bar{Q} |^2+\|K-\bar{K}\|^2), \\
&&|g( t,X,P,Y,Q,K) -g( t,\bar{X},\bar{P},\bar{Y},\bar{Q},\bar{K})|^{2}\\
&\leq& c(|X-\bar{X}|^{2}+|P-\bar{P}|^{2}+|Q-\bar{Q}|^{2}+\|K-\bar{K}\|^2)
+\gamma|Y-\bar{Y}|^2, \\
&&|G( t,X,P,Y,Q,K) -G( t,\bar{X},\bar{P},\bar{Y},\bar{Q},\bar{K})|
^{2}\\
&\leq& c( |X-\bar{X}| ^{2}+|P-\bar{P}|^{2}+|Y-\bar{Y}|^{2})
+\gamma(|Q-\bar{Q}|^{2}+\|K-\bar{K}\|^{2}),\\
&&|h( t,X,P,Y,Q,K) -h( t,\bar{X},\bar{P},\bar{Y},\bar{Q},\bar{K}) |
^{2}\\
&\leq& c( |X-\bar{X}| ^{2}+|P-\bar{P}|^{2}+|Y-\bar{Y}|^{2}
+|Q-\bar{Q} |^2+\|K-\bar{K}\|^2), \\
&&| \Psi (P) -\Psi ( \bar{P}) | \leq c|P-\bar{P}| ,\quad | \Phi (X)
-\Phi ( \bar{X}) | \leq c|X-\bar{X}|.
\end{eqnarray*}
\end{enumerate}
Then we claim the main result of this section.
\begin{theorem}\label{thm:3.2}
Under the assumptions (H2)$-$(H5), (\ref{eq:2}) has a unique
adapted solution $(X_t,P_t$, $Y_t,Q_t,K_t)\in M^{2}( 0,T;$
$\mathbb{R}^{n+m+n\times l+m\times d})\times
F^2_N(0,T;\mathbb{R}^m)$.
\end{theorem}
The proof of this theorem is divided into two parts, i.e., existence
and uniqueness. We first give the proof of uniqueness.

\noindent{\bf Proof.}\ (Uniqueness)\ Let $U_s=(X_s,P_s,Y_s,Q_s,K_s)$ and
$U'_s=(X'_s,P'_s,Y'_s,Q'_s,K'_s)$ be two solutions of (\ref{eq:2}). We set
$\widehat{U}=(\widehat{X},\widehat{P},
\widehat{Y},\widehat{Q},\widehat{K})$ $=(X-X',P-P',Y-Y',Q-Q',K-K')$.
Applying It\^{o}'s formula to $\langle
H\widehat{X}_s,\widehat{P}_s\rangle$ on $[0,T]$, we have
\begin{eqnarray*}
 &&E\langle H\widehat{X}_T,\Phi (X_T) -\Phi
(X_T^{\prime }) \rangle -E\langle H( \Psi ( P_0) -\Psi ( P_0^{\prime
}) ) ,\widehat{P}_0\rangle\\
&=&E\int_0^T\langle A(s,U_s)-A(s,U'_s),\widehat{U}_s\rangle ds\\
&\leq &-\mu _{1}E \int_{0}^{T}(| H\widehat{X}_s|
^{2}+|H\widehat{Y}_s|^{2})ds\\
&&-\mu _{2} E\int_0^T(|H^{\rm
T}\widehat{P}_s| ^{2}+|H^{\rm T}\widehat{Q}_s|^{2}+\|H^{\rm T}
\widehat{K}_s\|^{2})ds.
\end{eqnarray*}
Then
\begin{eqnarray*}
&&\mu _{1}E \int_{0}^{T}(| H\widehat{X}_s|
^{2}+|H\widehat{Y}_s|^{2})ds\\
&&+\mu _{2} E\int_0^T(|H^{\rm
T}\widehat{P}_s| ^{2}+|H^{\rm T}\widehat{Q}_s|^{2}+\|H^{\rm T}
\widehat{K}_s\|^{2})ds\leq 0.
\end{eqnarray*}

If $m>n$ and $\mu _{1}>0$, then we have $|H \widehat{X}_t |
^{2}\equiv 0$, $| H\widehat{Y}_t | ^{2}\equiv 0.$  Thus $X_t\equiv
X'_t, Y_t\equiv Y'_t.$ In particular, $\Phi(X_T)=\Phi(X'_T).$
Consequently, from Proposition \ref{pro:2.2}, it follows that $P_t\equiv P'_t,
Q_t\equiv Q'_t$ and $K_t\equiv K'_t.$

If $m<n$ and  $\mu _{2}>0$, then we have  $P_t\equiv P'_t, Q_t\equiv
Q'_t,$ and $K_t\equiv K'_t.$ In particular, $\Psi(P_0)=\Psi(P'_0).$
Thus from Proposition \ref{pro:2.2}, it follows that $X_t\equiv X'_t,$ and
$Y_t\equiv Y'_t.$

Similarly to above arguments, the desired result can be obtained
easily in the case $m=n$. \quad $\Box$
\begin{remark}\label{rmk:3.3}
In the proof of the uniqueness and existence, (H2) and (H3) can be
replaced by
\begin{enumerate}
\item[(H2)$'$]$\forall U=(X,P,Y,Q,K)$,
$\bar{U}=(\bar{X},\bar{P},\bar{Y},\bar{Q},\bar{K})\in
\mathbb{R}^{n+m+n\times l+m\times d}\times L^2_{\lambda(\cdot)}(
\mathbb{R}^m)$, $\forall t\in \left[ 0,T\right]$,
\begin{eqnarray*}
&&\langle A( t,U ) -A( t,\bar{U}) ,U -\bar{U}\rangle\\
&\geq& \mu _{1}(|H(X-\bar{X})|^{2}+| H(Y-\bar{Y})|^{2})\\
&&+\mu_{2}(|H^{\rm T}(P-\bar{P})|^{2}+|H^{T}(Q-\bar{Q})|^{2}+\|H^{\rm
T}(K-\bar{K})\|^{2}),
\end{eqnarray*}
\item[(H3)$'$] $\langle \Psi( P) -\Psi( \bar{P}), H^{\rm T}( P-\bar{P})
\rangle \geq \beta _{2}| H^{\rm
T}(P-\bar{P})|^{2},\quad \forall P,\bar{P}\in \mathbb{R}^{m}, \\
\langle \Phi (X) -\Phi ( \bar{X}), H(X-\bar{X}) \rangle \leq -\beta
_{1}| H(X-\bar{X})| ^{2},\quad \forall X,\bar{X}\in \mathbb{R}^{n}.$
\end{enumerate}
Here, $\mu _{1}$, $\mu _{2}$, $\beta _{1}$ and $\beta _{2}$ are
given nonnegative constants with $\mu _{1}+\mu _{2}>0$, $\beta
_{1}+\beta _{2}>0$, $\mu_{1}+\beta _{2}>0$, $\mu _{2}+\beta _{1}>0.$
Moreover we have $\mu _{1}>0$, $\beta _{1}>0$ (resp. $\mu _{2}>0$,
$\beta _{2}>0$ ) when $m>n$ (resp. $n>m$).
\end{remark}

The proof of the existence is a combination of the above technique
and the method of continuation systemically introduced by Yong \cite{Y}.
We divide the proof of existence into three cases: $m>n$, $m<n$ and
$m=n$.

{\bf Case 1}\quad  $m>n$

If $m>n$, then $\mu _1>0$ and $\beta _1>0$. We consider the
following family of FBDSDEP parametrized by $\alpha \in [0,1]$:
\begin{equation}\label{eq:3}
\left\{
\begin{array}{l}
dX_t =\left[ \alpha f(t,U_t) +f_0( t) \right] dt +\left[
\alpha g( t,U_t)+g_0(t) \right]dW_t-Y_tdB_t\\
\qquad\quad  +\int_Z\left[ \alpha h(t_-,U_{t_-},z)
+h_0( t_-,z)\right]\widetilde{N}(dzdt),   \\
dP_t =\left[ \alpha F( t,U_t) -( 1-\alpha
) \mu_1HX_t+F_0( t) \right]dt+Q_tdW_t  \\
\qquad\quad  +\left[ \alpha G( t,U_t) -( 1-\alpha ) \mu_1HY_t+G_0(
t) \right]dB_t
+\int_ZK_{t_-}(z)\widetilde{N}(dzdt),   \\
X_0 =\alpha \Psi (P_0)+\psi ,\quad P_T =\alpha \Phi (X_T) +(
1-\alpha ) HX_T+\varphi ,
\end{array}\right.
\end{equation}
where $U_t=(X_t,P_t,Y_t,Q_t)$, $( F_0,f_0,G_0,g_0,h_0) \in M^2(
0,T;\mathbb{R}^{m+n+m\times l+n\times d})\times F^2_N
(0,T;\mathbb{R}^n)$, $\psi \in L^2( \Omega ,{\cal
F}_0,P;\mathbb{R}^n)$ and $\varphi \in L^2( \Omega ,{\cal
F}_T,P;\mathbb{R}^m)$ are given arbitrarily.

When $\alpha =1$ the existence of the solution of (\ref{eq:3}) implies
clearly that of (\ref{eq:2}). Due to Proposition \ref{pro:2.2}, when $\alpha =0$, the
equation (\ref{eq:3}) is uniquely solvable. The following apriori lemma is
a key step in the proof of the method of continuation. It shows that
for a fixed $\alpha =\alpha _0\in [0,1),$ if (\ref{eq:3}) is uniquely
solvable, then it is also uniquely solvable for any $\alpha
\in[\alpha_0,\alpha_0+\delta_0]$, for some positive constant
$\delta_0$ independent of $\alpha_0.$

\begin{lemma}\label{lem:3.4}
We assume $m>n$. Under the assumptions (H2)-(H5), there exists
a positive constant $\delta_0$ such that if, apriori, for each
$\psi\in L^2( \Omega ,{\cal F}_0,P;\mathbb{R}^n) $, $\varphi \in
L^2( \Omega ,{\cal F} _T,P;\mathbb{R}^m) $, and $(
F_0,f_0,G_0,g_0,h_0) \in M^2( 0,T;\mathbb{R}^{m+n+m\times l+n\times
d})\times F^2_N(0,T$; $\mathbb{R}^m)$, (\ref{eq:3}) is uniquely solvable
for some $\alpha _0\in[0,1)$, then for each $\alpha \in [\alpha
_0,\alpha _0+\delta _0]$, and $ \psi \in L^2( \Omega ,{\cal
F}_0,P;\mathbb{R}^n) $, $\varphi \in L^2( \Omega ,{\cal
F}_T,P;\mathbb{R}^m) $, $( F_0,f_0,G_0,g_0,h_0) \in M^2(
0,T$; $\mathbb{R}^{m+n+m\times l+n\times d})\times
F^2_N(0,T;\mathbb{R}^m)$, (\ref{eq:3}) is also uniquely solvable in
$M^2( 0,T$;\\ $\mathbb{R}^{n+m+n\times l+m\times d})\times
F^2_N(0,T;\mathbb{R}^m)$.
\end{lemma}

{\bf Proof.}\ Since, for each $\psi \in L^2( \Omega ,{\cal
F}_0,P;\mathbb{R}^n)$, $\varphi \in L^2( \Omega ,{\cal
F}_T,P;\mathbb{R}^m) $, $( F_0,f_0$, $G_0,g_0,h_0)\in M^2(
0,T;\mathbb{R}^{m+n+m\times l+n\times d})\times
F^2_N(0,T;\mathbb{R}^m)$ there exists a unique solution of (\ref{eq:3}) for
$\alpha =\alpha _0$. Thus, for each $\bar U=( \bar X,\bar P,\bar Y,
\bar Q, \bar K)\in M^2( 0,T$; $\mathbb{R}^{n+m+n\times l+m\times
d})$ $\times F^2_N(0,T;\mathbb{R}^m) $, there exists a unique
$U=(X,P,Y,Q,K) \in M^2( 0,T;\mathbb{R}^{n+m+n\times l+m\times
d})\times F^2_N(0,T;\mathbb{R}^m)$ satisfying the following FBDSDEP:
\begin{eqnarray*}
dX_t &=&\left[ \alpha _0f( t,U_t) +\delta f( t, \bar U_t)
+f_0( t) \right]dt+\left[ \alpha _0g( t,U_t) +\delta g( t,\bar
U_t )+g_0( t) \right]d W_t\\
&&-Y_tdB_t+\int_Z\left[ \alpha _0 h(t_-,U_{t_-},z) +\delta h(t_-,\bar
U_{t_-},z) +h_0(
t_-,z)\right]\widetilde{N}(dzdt),\\
dP_t &=&\left[\alpha _0F( t,U_t) -( 1-\alpha _0) \mu
_1HX_t+\delta ( F( t,\bar U_t) +\mu _1H\bar X_t)
+F_0( t) \right]dt\\
&&+\left[ \alpha _0G( t,U_t) -(1-\alpha _0)\mu _1HY_t+\delta ( G(
t,\bar U_t) +\mu _1H\bar Y_t) +G_0( t)
\right]dB_t\\
&&+Q_tdW_t+\int_ZK_{t_-}(z)\widetilde{N}(dzdt),\\
X_0 &=&\alpha _0\Psi (P_0)+\delta \Psi (\bar P_0)+\psi,\\
P_T&=&\alpha _0\Phi (X_T) +( 1-\alpha _0) HX_T+\delta ( \Phi ( \bar X_T)
-H\bar X_T) +\varphi ,
\end{eqnarray*}
where $\delta \in ( 0,1)$ is independent of $\alpha _0$. We will
prove that the mapping defined by
\begin{eqnarray*}
&&(U_t,X_T,P_0)=I_{\alpha_0+\delta }(\bar U_t ,\bar X_T,\bar P_0):\\
&&M^2(0,T;\mathbb{R}^{n+m+n\times l+m\times d})
\times F^2_N(0,T;\mathbb{R}^m)\\
&&\times L^2( \Omega ,{\cal F}
_0,P;\mathbb{R}^n)\times L^2( \Omega ,{\cal F} _T,P;\mathbb{R}^m)\\
&&\to M^2( 0,T;\mathbb{R}^{n+m+n\times l+m\times d})
\times F^2_N(0,T;\mathbb{R}^m)\\
&&\times L^2( \Omega,{\cal F}_0,P;
\mathbb{R}^n)\times L^2( \Omega ,{\cal F} _T,P;\mathbb{R}^m ).
\end{eqnarray*}
is contractive for $\delta >0$ which is small enough.

Let
\begin{eqnarray*}
&&\bar U^{\prime}_t=(\bar X^{\prime }_t,\bar P^{\prime }_t,\bar
Q^{\prime }_t,\bar K^{\prime}_t) \in M^2(
0,T;\mathbb{R}^{n+m+n\times
l+m\times d})\times F^2_N(0,T;\mathbb{R}^m),\\
&&( X^{\prime }_t,P^{\prime }_t,Y^{\prime
}_t,Q^{\prime}_t,K^{\prime}_t) =U^{\prime}_t=I_{\alpha
_0+\delta }( \bar U^{\prime}_t),\\
&&\widehat{\bar U}_t =\bar U_t-\bar U^{\prime }_t=( \widehat{\bar
X}_t,\widehat{ \bar P}_t,\widehat{\bar Y}_t,\widehat{\bar
Q}_t,\widehat{\bar K}_t )\\
&&=( \bar X _t-\bar X^{\prime }_t,\bar
P_t-\bar P^{\prime }_t,\bar Y_t-\bar
Y^{\prime }_t,\bar Q_t-\bar Q^{\prime }_t, \bar K_t-\bar K^{\prime }_t ) , \\
&&\widehat{U}_t =U_t-U^{\prime }_t=(
\widehat{X}_t,\widehat{P}_t,\widehat{Y}_t, \widehat{Q}_t,
\widehat{K}_t)\\
&& =( X_t-X^{\prime }_t,P_t-P^{\prime }_t,Y_t-Y^{\prime
}_t ,Q_t-Q^{\prime }_t ,K_t- K^{\prime }_t).
\end{eqnarray*}
Applying It\^o's formula to $\langle H\widehat{X},\widehat{P}
\rangle$ on $[0,T]$ yields
\begin{eqnarray*}
&&E\langle H\widehat{X}_T,\alpha _0\widehat{\Phi }( X_T) +( 1-\alpha
_0) H\widehat{X}_T\rangle\\
&&-E \int_0^T\langle \alpha _0( A(s,U_s)
-A(s,U_s^{\prime})) ,\widehat U_s\rangle ds\\
&&+( 1-\alpha _0) \mu _1E\int_0^T( |
H\widehat{X}
_s| ^2+| H\widehat{Y}_s| ^2)ds \\
 &=&E\langle H\widehat{X}_T,\delta
H\widehat{\bar X} _T\rangle -E\langle H\widehat{X}_T,\delta
\widehat{\Phi } (\bar X_T)\rangle +E\langle H( \alpha
_0\widehat{\Psi }
(P_0)+\delta \widehat{\Psi }(\bar P_0)) ,\widehat{P}_0\rangle \\
&&+\delta E\int_0^T( \langle \widehat{P}_s,H\widehat{f} \rangle
+\langle H\widehat{X}_s,\widehat{F} \rangle +\langle
H\widehat{Y}_s,\widehat{G} \rangle+\langle
\widehat{Q}_s,H\widehat{g}\rangle+\langle \langle
\widehat{K}_s, H\widehat{h} \rangle \rangle )ds \\
&&+\delta \mu _1E\int_0^T( \langle H\widehat{X}_s,H\widehat{ \bar
X}_t\rangle +\langle H\widehat{Y}_s,H\widehat{\bar Y}_s\rangle )
ds,
\end{eqnarray*}
where
\begin{eqnarray*}
&&\widehat{f}= f( t,\bar U_t) -f( t,\bar U_t^{\prime
 }),\quad \widehat{g} =g( t,\bar U_t, \bar
K_t) -g( t,\bar U_t^{\prime }) ,\\
&&\widehat{F}=F( t,\bar U_t) -F( t,\bar U_t^{\prime }),\quad
\widehat{G} =G(
t,\bar U_t) -G( t,\bar U_t^{\prime }),\\
&&\widehat{h}=h( t,\bar U_t,\cdot) -h( t,\bar
U_t^{\prime},\cdot),\\
&&\widehat{\Psi}( \bar P_0)=\Psi( \bar P_0)-\Psi( \bar P^{\prime
}_0),\quad \widehat{\Psi}( P_0)=\Psi(P_0)-\Psi(
P^{\prime }_0),\\
&&\widehat{\Phi}( \bar X_T)=\Phi( \bar X_T)-\Phi( \bar X^{\prime
}_T),\quad\widehat{\Phi}(X_T)=\Phi( X_T)-\Phi( X^{\prime }_T).
\end{eqnarray*}
Let $B^2=M^2(0,T;\mathbb{R}^{n+m+n\times l+m\times d})\times
F^2_N(0,T;\mathbb{R}^m)$. Noting $m>n$, by virtue of (H2)-(H5), we
can easily deduce
\begin{eqnarray}\label{eq:4}
\nonumber&&( 1-\alpha _0+\alpha _0\beta _1) E| H\widehat{X}_T| ^2
+\mu _1E\int_0^T( |
H\widehat{X}_s| ^2+| H\widehat{Y}_s| ^2)ds\\
\nonumber&\leq& \delta CE\int_0^T( \| \widehat{U}_s\|_{B^2}^2+\|
\widehat{\bar U}_s\|_{B^2} ^2)ds\\
&&+\delta C( E{\bf \mid }\hat
P_0{\bf \mid }^2+E{\bf \mid }\widehat{X}_T{\bf \mid }^2+E{\bf \mid
}\widehat{\bar P}_0{\bf \mid }^2+E{\bf \mid }\widehat{\bar X}_T {\bf
\mid }^2).
\end{eqnarray}
 with some constant $C>0$. Hereafter, $C$ will be some
generic constant, which can be different from line to line and
depends only on the Lipschitz constants $c$, $\gamma$, $\mu _1$,
$\beta _1$, $H$ and $T$. It is obvious that $1-\alpha _0+\alpha
_0\beta _1\geq \beta$, $\beta=\min(1,\beta_1)>0.$

On the other hand, for the difference of the solutions $(
\widehat{P},\widehat{Q},\widehat{K}) =( P-P^{\prime}, Q-Q^{\prime},
K-K^{\prime })$, we apply a standard method of estimation. Applying
Proposition \ref{pro:2.3} to $| \widehat{P}_t| ^2$ on $[ t,T]$,  we get
\begin{eqnarray*}
&&E| \widehat{P}_t| ^2+E\int_t^T| \widehat
Q_s| ^2ds+E\int_t^T\|\widehat K_s\|^2ds \\
&=&E| \alpha _0\widehat{\Phi }(X_T) +(1-\alpha _0)H
\widehat{X}_T+\delta ( \widehat{\Phi }( \bar X_T) -H\widehat{
\bar X}_T) | ^2 \\
&&-2E\int_t^T\langle \widehat{P}_s,\alpha _0\widehat{F}( s,U_s) -(
1-\alpha _0) \mu _1H\widehat{X}_s+\delta ( \widehat{F}( s,\bar U_s)
+\mu _1H\widehat{\bar X}_s)\rangle ds \\
&&+E\int_t^T| \alpha _0\widehat{G}( s,U_s) -(1-\alpha _0)\mu
_1H\widehat{Y}_s+\delta ( \widehat{G}( s,\bar U_s) +\mu
_1H\widehat{\bar Y}_s) | ^2ds.
\end{eqnarray*}
By virtue of (H5), we have
\begin{eqnarray*}
&&E|\widehat{P}_t|^2+\frac{1-\gamma }4E\int_t^T( |\widehat Q_s|
^2+\| \widehat{K}_s\| ^2 )ds\\
&\leq& CE\int_t^T| \widehat P_s| ^2ds+C( E| \widehat{X}_T|
^2+\delta E| \widehat{\bar X}_T| ^2)\\
&&+CE\int_t^T( |\widehat{ X}_s|
^2+| \widehat{Y}_s| ^2+\delta \| \widehat{\bar U}_s\|_{B^2} ^2)ds.
\end{eqnarray*}
Then we can deduce
\begin{eqnarray}\label{eq:5}
\nonumber&&E| \widehat{P}_0| ^2+E\int_0^T( | \widehat{P}_s| ^2+|\widehat Q_s|
^2+\| \widehat{K}_s\| ^2 )ds\\
&\leq& C( E| \widehat{X}_T|
^2+\delta E| \widehat{\bar X}_T| ^2) +CE\int_0^T( |\widehat{ X}_s|
^2+| \widehat{Y}_s| ^2+\delta \| \widehat{\bar U}_s\|_{B^2} ^2)ds.
\end{eqnarray}
Combining the above two estimates (\ref{eq:4}) and (\ref{eq:5}), for a
sufficiently large constant $C>0$ we can easily have
\begin{eqnarray*}
&&E\int_0^T\| \widehat{U}_s\|_{B^2}^2ds+E| \widehat{X} _T| ^2+E|
\widehat{P}_0| ^2\\
&\leq&\delta C(E \int_0^T\|\widehat{\bar
U}_s\|_{B^2}^2ds+E| \widehat{\bar X} _T| ^2+E| \widehat{\bar
P}_0| ^2),
\end{eqnarray*}
We now choose $\delta _0=\displaystyle\frac 1{2C}$. It is clear
that, for each fixed $\delta \in \left[ 0,\delta _0\right] $, the
mapping $I_{\alpha_0+\delta }$ is contract in the sense that
\begin{eqnarray*}
&&E\int_0^T\| \widehat{U}_s\|_{B^2}^2ds+E| \widehat{X} _T| ^2+E|
\widehat{P}_0| ^2\\
&\leq& \displaystyle\frac 12(E
\int_0^T\|\widehat{\bar U}_s\|_{B^2}^2ds+E| \widehat{\bar X}
_T| ^2+E| \widehat{\bar P}_0| ^2),
\end{eqnarray*}
Thus, this mapping has a unique fixed point $U=(X,P,Y,Q,K) \in
M^2(0,T$; $\mathbb{R}^{n+m+n\times l+m\times d})$ $\times
F_N^2(0,T;\mathbb{R}^m) $, which is the solution  of (\ref{eq:3}) for
$\alpha =\alpha_0+\delta $, $\delta\in [0,\delta_0]$.  \quad $\Box$

{\bf Case 2}\quad $m<n$

If $m<n$, then $\mu_{2}>0$ and $\beta _{2}>0$. We consider following
equations:
\begin{equation}\label{eq:6}
\left\{
\begin{array}{l}
dX_t =\left[ \alpha f( t,U_t) -( 1-\alpha ) \mu _2H^{\rm
T}P_t+f_0( t) \right]dt\\
\qquad\quad +\left[ \alpha g( t,U_t) -( 1-\alpha )
\mu _2H^{\rm T}Q_t+g_0( t) \right] dW_t-Y_tdB_t\\
\qquad\quad +\displaystyle\int_Z\left[ \alpha
h(t_-,U_{t_-},z) +( 1-\alpha) \mu
_2H^{\rm T}K_{t_-}(z)+h_0(t_-,z)\right]\widetilde{N}(dzdt),  \\
dP_t =\left[ \alpha F( t,U_t) +F_0( t) \right] dt+Q_tdW_t+\left[ \alpha G( t,U_t) +G_0( t)
\right]dB_t\\
\qquad\quad -\displaystyle\int_ZK_{t_-}(z)\widetilde{N}(dzdt),   \\
X_0 =\alpha \Psi ( P_0) +( 1-\alpha ) H^{\rm T}P_0+\psi ,\quad P_T
=\alpha \Phi ( X_T) +\varphi .
\end{array}
\right.
\end{equation}

Due to Proposition \ref{pro:2.2}, when $\alpha =0$, the equation (\ref{eq:6}) is
uniquely solvable. When $\alpha =1$ the existence of the solution of
(\ref{eq:6}) implies clearly that of (\ref{eq:2}). By the techniques similar to
Lemma \ref{lem:3.4}, We can also prove the following lemma.

\begin{lemma}\label{lem:3.5}
We assume $m<n$. Under assumptions (H2)-(H5), there exists a
positive constant $\delta_0$ such that if, apriori, for each
$\psi\in L^2( \Omega ,{\cal F}_0,P;R^n) $, $\varphi \in L^2( \Omega
,{\cal F} _T,P;R^m) $, and $( F_0,f_0,G_0,g_0,h_0) \in M^2(
0,T;R^{m+n+m\times l+n\times d})\times F^2_N(0,T$; $\mathbb{R}^m)$,
(\ref{eq:6}) is uniquely solvable for some $\alpha _0\in[0,1)$, then
for each $\alpha \in [\alpha _0,\alpha _0+\delta _0]$, and $ \psi
\in L^2( \Omega ,{\cal F}_0,P;R^n) $, $\varphi \in L^2( \Omega
,{\cal F}_T,P;R^m) $, $( F_0,f_0,G_0,g_0,h_0) \in M^2(
0,T$; $R^{m+n+m\times l+n\times d})\times F^2_N(0,T;\mathbb{R}^m)$,
(\ref{eq:6})is also uniquely solvable in $M^2( 0,T$;\\
$R^{n+m+n\times l+m\times d})\times F^2_N(0,T;\mathbb{R}^m)$.
\end{lemma}

{\bf Case 3}\hspace{2mm}  $m=n$

From (H2) and (H3), we note that we only need to consider two cases
as follows:

(1)\ If $\mu _1>0$, $\mu _2\geq 0$, $\beta _1>0$, and $\beta _2\geq
0$, we can have the same result as Lemma \ref{lem:3.4}.

(2)\ If $\mu _1\geq 0$, $\mu _2>0$, $\beta _1\geq 0$, and $\beta
_2>0$, the same result as Lemma \ref{lem:3.5} holds.

\vskip 1mm Now we can give proof of the existence of Theorem \ref{thm:3.2}.

{\bf Proof of the existence of Theorem \ref{thm:3.2}.}\
For Case 1, we know that, for each $\psi \in L^2( \Omega ,{\cal
F}_0,P;R^n)$, $ \varphi \in L^2( \Omega ,{\cal F}_T,P;R^m) $, $(
F_0,f_0,G_0,g_0,h_0) \in M^2( 0,T$; $R^{m+n+m\times l+n\times d})\times
F^2_N(0,T;\mathbb{R}^m)$, (\ref{eq:3}) has a unique solution as $\alpha
=0$. It follows from Lemma \ref{lem:3.4} that there exists a positive constant
$\delta _0= \delta _0(c,\gamma,\beta_1,\mu _1$, $H,T)$ such that for
any $\delta \in [0,\delta _0]$ and $\psi \in L^2( \Omega ,{\cal
F}_0,P;R^n)$, $\varphi \in L^2( \Omega ,{\cal F}_T,P$; $R^m) $, $(
F_0,f_0,G_0,g_0,h_0) \in M^2( 0,T;R^{m+n+m\times l+n\times d})\times
F^2_N(0,T;\mathbb{R}^m)$, (\ref{eq:3}) has a unique solution for $\alpha
=\delta $. Since $\delta _0$ depends only on $c$, $\gamma$,
$\beta_1$, $\mu _1$, $H$, and $T $, we can repeat this process for
$N$ times with $1\leq N\delta _0<1+\delta _0$. In particular, for
$\alpha =1$ with $( F_0,f_0,G_0,g_0,h_0) \equiv 0$, $\varphi \equiv
0$, $\psi \equiv 0$, (\ref{eq:3}) has a unique solution in $M^2(
0,T;R^{n+m+n\times l+m\times d})$ $\times F^2_N(0,T;\mathbb{R}^m)$.

For Case 2, we know that, for each $\psi \in L^2( \Omega ,{\cal
F}_0,P;\mathbb{R}^n) $, $ \varphi \in L^2( \Omega ,{\cal F}_T$,
$P$; $\mathbb{R}^m) $, $(
F_0,f_0,G_0,g_0,h_0) \in M^2( 0,T$; $R^{m+n+m\times l+n\times
d})\times F^2_N(0,T;\mathbb{R}^m)$, (\ref{eq:6}) has a unique solution as
$\alpha =0$. It follows from Lemma \ref{lem:3.5} that there exists a positive
constant $\delta _0= \delta _0(c,\gamma,\beta_2,\mu_2,H,T)$ such
that for any $\delta \in [0,\delta _0]$ and $\psi \in L^2( \Omega
,{\cal F}_0,P;R^n)$, $\varphi \in L^2( \Omega ,{\cal F}_T,P;R^m) $,
$( F_0,f_0,G_0,g_0,h_0) \in M^2( 0,T$; $R^{m+n+m\times l+n\times
d})\times F^2_N(0,T$; $\mathbb{R}^m)$, (\ref{eq:6}) has a unique solution for
$\alpha =\delta $. Since $\delta _0$ depends only on $c$, $\gamma$,
$\beta_2$, $\mu _2$, $H$, and $T $, we can repeat this process for
$N$ times with $1\leq N\delta _0<1+\delta _0$. In particular, for
$\alpha =1$ with $( F_0,f_0,G_0,g_0,h_0) \equiv 0$, $\varphi \equiv
0$, $\psi \equiv 0$, (\ref{eq:6}) has a unique solution in $M^2(
0,T;R^{n+m+n\times l+m\times d})$ $\times F^2_N(0,T;\mathbb{R}^m)$.

Similar to the above two cases, the desired result can be obtained
in Case 3.  \quad $\Box$

\section{FBDSDEP depending on parameters}
\label{sec:4}

In this section, the continuity of the solutions to FBDSDEP
depending on parameters is discussed. Let $\{ f_\alpha ,g_\alpha$,
$h_\alpha ,F_\alpha ,G_\alpha ,\Psi _\alpha ,\Phi _\alpha, \alpha
\in R \}$ be a family of data of FBDSDEP as follows
\begin{equation}\label{eq:7}
\left\{
\begin{array}{l}
dX_t^\alpha =f_\alpha ( t,X_t^\alpha, P_t^\alpha, Y_t^\alpha,
Q_t^\alpha,  K_t^\alpha )dt+g_\alpha (
 t,X_t^\alpha, P_t^\alpha,
Y_t^\alpha, Q_t^\alpha,  K_t^\alpha) dW_t\\
\qquad\quad  -Y_t^\alpha dB_t+\displaystyle\int_Z h_\alpha ( t_-,X_{t_-}^\alpha,
P_{t_-}^\alpha,
Y_{t_-}^\alpha, Q_{t_-}^\alpha, K_{t_-}^\alpha(z),z ) \widetilde{N}(dzdt), \\
dP_t^\alpha =F_\alpha (  t,X_t^\alpha, P_t^\alpha, Y_t^\alpha,
Q_t^\alpha,  K_t^\alpha)dt+G_\alpha (
 t,X_t^\alpha, P_t^\alpha,
Y_t^\alpha, Q_t^\alpha,  K_t^\alpha ) dB_t\\
\qquad\quad  +Q_t^\alpha W_t
+\displaystyle\int_ZK_{t_-}(z)\widetilde{N}(dzdt), \\
X_0 =\Psi_\alpha  ( P_0),\quad P_T =\Phi_\alpha  ( X_T),
\end{array}\right.
\end{equation}
with solutions denoted by $(X^\alpha, P^\alpha, Y^\alpha, Q^\alpha,
K^\alpha ).$

 Let us give some assumptions:
\begin{enumerate}
\item[(H6)]
\begin{enumerate}
\item[(i)] The family $\left\{ f_\alpha, F_\alpha, g_\alpha ,G_\alpha
,h_\alpha,\Psi _\alpha ,\Phi _\alpha ;\alpha \in R\right\}$ satisfy
the equi-\\
Lipschitz conditions with the same constant $c$ as in (H5);
\item[(ii)] $\left\{ f_\alpha, F_\alpha, g_\alpha, G_\alpha ,h_\alpha
,\Psi _\alpha ,\Phi _\alpha \right\}$ are continuous for $\alpha$,
in their existing\\
space norm sense respectively.
\end{enumerate}
\end{enumerate}

Then we have the following continuity result.

\begin{theorem}\label{thm:4.1}
 Let  $\left\{ f_\alpha, F_\alpha ,g_\alpha, G_\alpha, h_\alpha,
 \Psi _\alpha ,\Phi _\alpha ;\alpha \in \mathbb{R}\right\}
$  be a family of FBDSDEP (\ref{eq:7}) satisfying (H2)$-$(H6)
with solutions denoted by $( X^\alpha,P^\alpha ,Y^\alpha ,Q^\alpha
,K^\alpha) $. Then, the family of functions $\{ ( X^\alpha,
P^\alpha, Y^\alpha, Q^\alpha, K^\alpha, X_T^\alpha ,P_0^\alpha )$;
$\alpha \in R\}$  is continuous for $\alpha $ in $M^2(0,T$;
$\mathbb{R}^{n+m+n\times l+m\times d})\times
F_N^2(0,T;\mathbb{R}^m)\times L^2( \Omega ,{\cal
F}_T,P;\mathbb{R}^n)
 \times L^2( \Omega,$ ${\cal F}_T,P;\mathbb{R}^m)$ .
\end{theorem}

{\bf Proof.}\ For notational convenience, we only prove the
continuity of FBDSDEP (\ref{eq:7}) at $\alpha =0$. We need to get that $(
X^\alpha, P^\alpha, Y^\alpha$, $Q^\alpha, K^\alpha, X_T^\alpha,
P_0^\alpha)$  converges to  $( X^0, P^0$, $Y^0, Q^0, K^0, X_T^0,
P_0^0)$ in $M^2( 0,T;\mathbb{R}^{n+m+n\times l+m\times d})
 \times F_N^2(0,T;\mathbb{R}^m)
 \times L^2( \Omega ,{\cal F}_T,P;\mathbb{R}^n)
 \times L^2( \Omega,{\cal F}_T,$ $P;\mathbb{R}^m)
 $ as $\alpha \rightarrow 0$.

We set
\begin{eqnarray*}
\widehat{U}_t &=& U_t^\alpha -U_t^0=( \widehat{X}_t,\widehat{P}_t,
\widehat{Y}_t,\widehat{Z}_t,\widehat{K}_t)\\
&=&( X_t^\alpha
-X_t^0,P_t^\alpha -P_t^0,Y_t^\alpha -Y_t^0,Q_t^\alpha -Q_t^0,
K_t^\alpha -K_t^0 ).
\end{eqnarray*}
Thus
\begin{eqnarray*}\left\{
\begin{array}{l}
d\widehat{X}_t=( f_\alpha ( t,U_t^\alpha ) -f_0( t,U_t^0) )
dt+( g_\alpha ( t,U_t^\alpha)-g_0( t,U_t^0) )
dW_t-\widehat{Y}_tdB_t\\
\qquad\quad+\displaystyle\int_Z( h_\alpha ( {t_-},U_{t_-}^\alpha,z)
-h_0( t_-,U_{t_-}^0z) )\widetilde{N}(dzdt),\\
d\widehat{P}_t =( F_\alpha ( t,U_t^\alpha ) -F_0( t,U_t^0) )
dt+( G_\alpha ( t,U_t^\alpha )-G_0( t,U_t^0) ) dB_t\\
\qquad\quad+\widehat{Q}_tdW_t
+\displaystyle\int_Z\widehat{K}_{t_-}^\alpha(z)\widetilde{N}(dzdt),\\
\widehat{X}_0 =\Psi _\alpha ( P_0^\alpha ) -\Psi _0( P_0^0) , \quad
\widehat{P}_T =\Phi _\alpha ( X_T^\alpha ) -\Phi _0( X_T^0).
\end{array}\right.
\end{eqnarray*}
 Applying It\^o's formula to $|\widehat{X}_t|^2$,
$|\widehat{P}_t|^2$, and $\left<
H\widehat{X}_t,\widehat{P}_t\right>$ with the usual technique
similarly to Lemma \ref{lem:3.4}, we can obtain
\begin{eqnarray*}
&&E\int_0^T\|\widehat{U}_s\|_{B^2}^2ds+E| \widehat{X} _T| ^2+E|
\widehat{P}_0| ^2\\
&\leq &CE\int_0^T( | f_\alpha (s,U_s^0) -f_0(
s,U_s^0) | ^2+| g_\alpha(s,U_s^0) -g_0( s,U_s^0)|^2)ds\\
&&+CE\int_0^T( | F_\alpha (s,U_s^0) -F_0(s,U_s^0) | ^2+|
G_\alpha(s,U_s^0) -G_0(s,U_s^0) | ^2)ds \\
&&+CE\int_0^T\|h_\alpha (s,U_s^0,\cdot) -h_0(s,U_s^0,\cdot)\|^2ds\\
&& +CE|\Phi_\alpha(y_T^0)-\Phi_0(y_T^0)|^2+CE| \Psi_\alpha (Y_0^0)
-\Psi_0(Y_0^0)|^2.
\end{eqnarray*}
where $C >0$ is some constant depending on the Lipschitz constants
$c$, $\gamma$, $\mu _1$, $\mu _2$, $\beta _1$, $\beta _2$, $H$, and
$T$. Therefore, from (H6), it follows that $(X^\alpha, P^\alpha,
Y^\alpha, Q^\alpha$, $K^\alpha, X_T^\alpha ,P_0^\alpha)$ converges
to $(X^0,P^0,Y^0,Q^0,K^0,X_T^0,P_0^0)$ in the space $M^2(
0,T$; $R^{n+m+n\times l+m\times d})\times F_N^2(0,T;\mathbb{R}^m)\times
L^2( \Omega,$ ${\cal F}_T,P;R^n) \times L^2( \Omega ,{\cal F}
_0,P;R^m)$ as $\alpha \rightarrow 0$.  \quad
$\Box$

In fact, we can also discuss the differentiability  of the solution
to FBDSDEP depending on parameters. The method is similar. These two
properties are important and make FBDSDEP be used widely especially
in practice.

\section{The probabilistic interpretation of SPDIEs}
\label{sec:5}

The connection of BDSDEs and systems of second-order quasilinear
SPDEs was observed by Pardoux and Peng \cite{PP2}. This can be regarded as
a stochastic version of the well-known Feynman-Kac formula which
gives a probabilistic interpretation for second-order SPDEs of
parabolic types. Thereafter this subject has attracted many
mathematicians, referred to Bally and Matoussi \cite{BM}, Zhang and Zhao
\cite{ZZ}, Hu and Ren \cite{HR}, see also Ren et al. \cite{RLH}]. In \cite{RLH}], the authors
got a probabilistic interpretation for the solution of a semilinear
SPDIE, via BDSDEs with L$\acute{\rm e}$vy process.  This section can
be viewed as a continuation of such a theme, and will exploit the
above theory of fully coupled FBDSDEP in order to provide a
probabilistic formula for the solution of a quasilinear SPDIE.

For each $x\in \mathbb{R}^n$, let $\left\{ X_t,P_t,Q_t,K_t;0\leq t<T
\right\}$ denote the solution of the FBDSDEP:
\begin{equation}\label{eq:8}
\left\{
\begin{array}{llll}
dX_t & = & f(t,X_t,P_t,Q_t,K_t)dt+g(t, X_t,P_t) dW_t+\displaystyle\int_Z h(t_-,X_{t_-},z)
\widetilde{N}(dzdt),\\
dP_t & = & F(t,X_t,P_t,Q_t,K_t) dt+G(t,X_t,P_t,Q_t,K_t)
dB_t\\
&&+Q_tdW_t+\displaystyle\int_Z K_{t_-}(z)\widetilde{N}(dzdt),\\
 X_t&=&x, \quad P_T=\Phi ( X_T),
\end{array}
\right.
\end{equation}
where
\begin{center}
$\begin{array}{lllll} F & :[0,T]\times \mathbb{R}^n\times
\mathbb{R}^m\times
 \mathbb{R}^{m\times d}\times L^2_{\lambda(\cdot)}(
\mathbb{R}^m) & \to & \mathbb{R}^m, \\
f & :[0,T]\times  \mathbb{R}^n\times\mathbb{R}^m\times
 \mathbb{R}^{m\times d}\times L^2_{\lambda(\cdot)}(
\mathbb{R}^m) & \to & \mathbb{R}^n, \\
G & :[0,T]\times  \mathbb{R}^n\times \mathbb{R}^m\times
 \mathbb{R}^{m\times d}\times L^2_{\lambda(\cdot)}(
\mathbb{R}^m) & \to & \mathbb{R}^{m\times l}, \\
g & :[0,T]\times  \mathbb{R}^n\times \mathbb{R}^m & \to &
\mathbb{R}^{n\times d},\\
h & :[0,T]\times  \mathbb{R}^n \times Z & \to & \mathbb{R}^{n},\\
\Phi & : \mathbb{R}^n  & \to & \mathbb{R}^m,
\end{array}
$
\end{center}
satisfy (H2)-(H5), and
\begin{enumerate}
\item[(H7)] $f$, $g$, $h$, $F$ and $G$  are of class $C^3$, and $\Phi$ is of class $C^2$.
\end{enumerate}

We now relate FBDSDEP (\ref{eq:8}) to the following system of quasilinear
second-order parabolic SPDIE:
\begin{eqnarray}\label{eq:9}
\left\{
\begin{array}{lllll}
{\cal L}u(t,x)dt\\
=F(t,x,u(t,x),\nabla u(t,x)g(t,x,u(t,x)),u(t,x+h(t,x,\cdot))-u(t,x))dt\\
+G(t,x,u(t,x),\nabla u(t,x)g(t,x,u(t,x)),u(t,x+h(t,x,\cdot))-u(t,x))dB_t,\\
(t,x)\in [0,T]\times R^n,\quad u(T,x)=\Phi (x),
\end{array}
\right.
\end{eqnarray}
where $u:\mathbb{R}_+\times\mathbb{R}^n\to \mathbb{R}^m$,
\begin{eqnarray*}
{\cal L}u =\left(\begin{array}{c} Lu_1\\\vdots\\
Lu_m\end{array}\right),
\end{eqnarray*}
with
\begin{eqnarray*}
&&Lu_k(t,x)\\
&=&\displaystyle\frac{\partial u_k}{\partial
t}(t,x)+\displaystyle\frac{1}{2}\sum\limits_{i,j=1}^n
(gg*)_{ij}(t,x,u(t,x))\displaystyle\frac{\partial^2 u_k }{\partial
x_i\partial x_j}(t,x)\\
&& +\sum\limits_{i=1}^nf_{i}(t,x,u(t,x),\nabla
u(t,x)g(t,x,u(t,x)),\\
&&u(t,x+h(t,x,\cdot))-u(t,x))
\displaystyle\frac{\partial u_k }{\partial x_i}(t,x)\\
&&+\displaystyle\int_Z(u_k(t,x+h(t,x,z))-u_k(t,x)-\sum\limits_{i=1}^n
h_i(t,x,z)\displaystyle\frac{\partial u_k}{\partial x_i}(t,x)
)\lambda(dz),\\
&&k=1,\cdots,m.
\end{eqnarray*}
We can assert that
\begin{theorem}\label{thm:5.1}
Assume that the functions $f$, $g$, $h$, $F$, $G$ and $\Phi$ in
FBDSDEP (\ref{eq:8}) are deterministic and that they satisfy the
assumptions (H2)-(H5) and (H7). Suppose SPDIE (\ref{eq:9})
has a unique solution $u(t,x)\in C^{1,2}(\Omega \times
[0,T]\times \mathbb{R}^n; \mathbb{R}^m)$. Then, for any given
$(t,x)$, $u(t,x)$ has the following interpretation
\begin{eqnarray}\label{eq:10}
u(t,x)= P_t,\end{eqnarray} where $P_t$ is determined uniquely by
FBDSDEP (\ref{eq:8}).
\end{theorem}
{\bf Proof.}Applying  Proposition \ref{pro:2.4} to $u(t,X_t)$, we obtain
\begin{eqnarray*}
&&u(T,X_T)-u(t,X_t)\\&=&\int_t^T\frac{\partial u}{\partial s}(s,X_s)ds+\int_t^T\sum\limits_{i=1}^nf_{i}(s,X_s,P_s,Q_s,K_s)
\displaystyle\frac{\partial u}{\partial x_i}(s,X_s)ds\\
&& +\int_t^T\nabla u(s,X_s)g(s,X_s,P_s)dW_s\\
&& +\int_t^T\displaystyle\int_Z(u(s_-,X_{s_-}
+h(s_-,X_{s_-},z))-u(s,X_{s_-}))\widetilde{N}(dzds)\\
&&+\int_t^T\displaystyle\frac{1}{2}\sum\limits_{i,j=1}^n
(gg*)_{ij}(s,X_s,P_s)\displaystyle\frac{\partial^2u }{\partial
x_i\partial x_j}(s,X_s)ds\\
&&+\int_t^T\displaystyle\int_Z(u(s,X_s+h(s,X_s,z))-u(s,X_s)\\
&&-\sum\limits_{i=1}^n
h_i(s,X_s,z)\displaystyle\frac{\partial u}{\partial x_i}
(s,X_s))\lambda(dz)ds.
\end{eqnarray*}
Let
\begin{eqnarray*}
&&(P_t,Q_t,K_t)\\
&=& (u(t,X_t),\nabla u(t,X_t)g(t,X_t,u(t,X_t)),
u(t,X_t+h(t,X_t,\cdot))-u(t,X_t)).
\end{eqnarray*}
 Because $u(t,x)$ satisfies SPDIE {\rm (\ref{eq:9})}, it holds that
\begin{eqnarray*}
&&\Phi(X_T)-u(t,X_t)\\
&=&\int_t^TF(s,X_s,u(s,X_s),\nabla
u(s,X_s)g(s,X_s,u(s,X_s)),\\
&&u(s,X_s+h(s,X_s,\cdot))-u(s,X_s))ds\\
&&+ \int_t^T G(s,X_s,u(s,X_s),\nabla
u(s,X_s)g(s,X_s,u(s,X_s)),\\
&&u(s,X_s+h(s,X_s,\cdot))-u(s,X_s))dB_s \\&& +\int_t^T\nabla u(s,X_s)g(s,X_s,u(s,X_s))dW_s\\
&&+\int_t^T\displaystyle\int_Z(u(s_-,X_{s_-}+h(s_-,X_{s_-},z))-u(s_-,X_{s_-}))\widetilde{N}(dzds).
\end{eqnarray*}
It is easy to check that $(u(t,X_t),\nabla
u(t,X_t)g(t,X_t,u(t,X_t)), u(t,X_t+h(t,X_t,\cdot))-u(t,X_t))$
coincides with the unique solution to BDSDEP of (\ref{eq:8}). It follows
that
$$u(t,x)= P_t.$$  \quad $\Box$

\begin{remark}\label{rmk:5.2}
(\ref{eq:10}) can be called a Feynman-Kac formula for SPDIE (\ref{eq:9}).
Furthermore, with regard to FBDSDEs driven by L$\acute{\rm e}$vy
process rather than by Poisson process, we can get the similar
Feynman-Kac formula as (\ref{eq:10}) for a SPDIE.
\end{remark}

\section{Example: a doubly stochastic Hamiltonian system with Brownian motions and
Poisson process}\label{sec:6}

Consider the following doubly stochastic Hamiltonian system with
Brownian motions and Poisson process
\begin{eqnarray}\label{eq:11}
\left\{
\begin{array}{llll}
dX_t & = & H_Pdt+H_QdW_t-Y_tdB_t+\displaystyle\int_Z H_K\widetilde{N}(dzdt),\\
dP_t & = & -H_Xdt-H_YdB_t+Q_tdW_t+\displaystyle\int_Z K_t(z)\widetilde{N}(dzdt),\\
 X_0&=&\Psi _P( P_0),\quad P_T=\Phi _X( X_T) ,
\end{array}
\right.
\end{eqnarray}
where $H(X,P,Y,Q,K):\mathbb{R}^4\times
L_{\lambda(\cdot)}^2(\mathbb{R})\rightarrow \mathbb{R}$, $\Phi ( X):
\mathbb{R}\rightarrow \mathbb{R}$, $\Psi ( P):\mathbb{R}\rightarrow
\mathbb{R}$; $H_P\doteq \nabla _PH,$ $\Phi _X\doteq \nabla _X\Phi
,\Psi _P\doteq \nabla _P\Psi $.
 The Brownian motions $\left\{ W_t;0\leq
t\leq T\right\} $ and $\left\{ B_t;0\leq t\leq T\right\} $ are both
assumed to be 1-dimensional. Assume that both the derivatives of
$2$-order of $H$ and the derivatives of $1$-order of $\Phi $ and
$\Psi $ are bounded, $H$ is concave on $( P,Q,K) $ and convex on
$(X,Y) $ in the following sense $\mu_1 >0$, and $\mu_2
>0$:
\begin{center}
$\left(
\begin{array}{ccccc}
-H_{XX} & -H_{XP} & -H_{XY} & -H_{XQ} & -H_{XK}\\
H_{PX} & H_{PP} & H_{PY} & H_{PQ} & H_{PK}\\
-H_{YX} & -H_{YP} & -H_{YY} & -H_{YQ} & -H_{YK}\\
H_{QX} & H_{QP} & H_{QY} & H_{QQ} & H_{QK}\\
H_{KX} & H_{KP} & H_{KY} & H_{KQ} & H_{KK}
\end{array}
\right) \leq -\left(
\begin{array}{ccccc}
\mu_1 & 0 & 0 & 0 &0\\
0 & \mu_2  & 0 & 0 & 0\\
0 & 0 &\mu_1  & 0 & 0\\
0 & 0 & 0 & \mu_2 & 0\\
0 & 0 & 0 & 0 & \mu_2
\end{array}
\right),$
\end{center}
where $\forall ( X,P,Y,Q,K) \in \mathbb{R}^4\times
L_{\lambda(\cdot)}^2(\mathbb{R})$, and  $\Phi $ is convex on $X$:
$\Phi _{XX}\geq 0,\ \forall X\in R$; $\Psi$ is concave on $P$: $\Psi
_{PP}\leq 0,\ \forall P\in R.$ By Theorem \ref{thm:3.2}, we claim that this
doubly stochastic Hamiltonian system (\ref{eq:11}) has a unique solution $(
X,P,Y,Q,K) $ in $M^2( 0,T;\mathbb{R}^4)\times
F^2_N(0,T;\mathbb{R})$.

\end{document}